\documentclass[12pt]{amsart} 
\usepackage{amsmath}
\usepackage{amssymb}

\newtheorem{theorem}{Theorem}[section]

\newtheorem{lemma}[theorem]{Lemma}
\newtheorem{corollary}[theorem]{Corollary}

\newtheorem{definitionn}[theorem]{Definition}

\def\C{{\mathbf C}}

\def\P{{\mathbf P}}

\def\x{{\mathbf x}}
\def\y{{\mathbf y}}

\def\cA{{\mathcal A}}
\def\cB{{\mathcal B}}
\def\cC{{\mathcal C}}

\def\cE{{\mathcal E}}

\def\cH{{\mathcal H}}

\def\cP{{\mathcal P}}

\def\cX{{\mathcal X}}

\def\ord{{\rm ord}}
\def\deg{{\rm deg}}

\def\min{{\rm min}}

\def\fq{{F_q}}

\begin{document}

\title[On Near-MDS Elliptic Codes]{On Near-MDS Elliptic Codes} 

\author[M. Giulietti]{Massimo Giulietti}\thanks{Massimo Giulietti is with the Dipartimento di Matematica e Informatica,
Universit\`a di Perugia, 06123 Perugia,
Italy. E-mail: giuliet@dipmat.unipg.it}\thanks{This research was performed within the activity of GNSAGA
of the Italian CNR, with the financial support of the Italian
Ministry MIUR, project `` Strutture geometriche, combinatoria e loro
applicazioni'', PRIN 2001-2002.}

\begin{abstract}
The Main Conjecture on maximum distance separable (MDS) codes
states that, except for some special cases, the maximum length of a $q$-ary linear MDS code of
is $q+1$. This conjecture
does not hold true for near
maximum distance separable codes because of the existence of
$q$-ary near-MDS elliptic codes
having length bigger than $q+1$. An
interesting related question
is whether a near-MDS
elliptic code may be extended to a longer near-MDS code. Our
results are some non-extendability results and an alternative and
simpler construction for certain known near-MDS elliptic codes.
\end{abstract}

\maketitle

{\bf Keywords:} Projective Spaces, Near-MDS Codes, Elliptic Curves.
\section{Introduction}
Let $\fq$ be a finite field with $q$ elements and $\fq^n$ the
vector space of $n$-tuples over $\fq$. A \underline{$q$-ary linear code}
$\C$ of length $n$ and dimension $k$ is a $k$-dimensional subspace
of $\fq^n$. The number of non-zero positions in a vector $\x\,\in
\C$ is called the \underline{Hamming weight} $w(\x)$ of $\x$; the \underline{Hamming distance} $d(\x,\y)$ between two vectors $\x,\y\,\in \,
\C$ is defined by $d(\x,\y)=w(\x - \y)$. The \underline{minimum
distance} of $\C$ is $$ d(\C):=\min\{w(\x)\mid \x \,\in\, \C, \,
\x \neq 0\}\,, $$ and a $q$-ary linear code of length $n$,
dimension $k$ and minimum distance $d$ is indicated as an $[n,k,d]_q$ code. For such
codes the Singleton bound holds: $$ d\le n-k+1\,. $$ The non-negative integer
$s(\C):=n-k+1-d$ is referred to as the \underline{Singleton defect} of $\C$.

A linear code $\C$ with $s(\C)=0$ is said to be \underline{maximum
distance separable}, or briefly MDS. A code with $s(\C)=1$ is
called \underline{almost-MDS}, or AMDS for short. The dual $\C^\perp$ of
a code $\C$ consists of all the vectors of $\fq^n$ orthogonal to
every codewords in $\C$:
$$
 \C^\perp:=\{\x\,\in\, \fq^n\mid
\langle\x,\y\rangle=0 \,\,\,{\rm for\,\, any}\,\, \y\,\in\,\C\}\,,
$$ where $\langle,\rangle$ denotes the inner product in $\fq^n$.
Unlike the MDS case, the dual of an AMDS code need not be AMDS.
This motivates to define $\C$ to be \underline{near-MDS} (NMDS) when
$s(\C)=s(\C^\perp)=1$.


For given $k$ and $q$, let $m(k,q)$ be the maximum length of a $q$-ary
linear MDS code of dimension $k$. The Main Conjecture on MDS codes states that
$m(k,q)=q+1$ provided that $2\le k < q$, except for the case $m(3,q)=m(q-1,q)=q+2$
for even $q$ (see e.g. \cite[p. 13]{Gtsf-vla}). The situation is quite different for NMDS codes, since $q$-ary linear
NMDS codes of length bigger than $q+1$ arise from elliptic curves via Goppa
construction. In particular the following theorem
holds (\cite[Sec. 3.2]{Gtsf-vla}).


        \begin{theorem}\label{Th1.1}
Let $q=p^m$, $p$ prime. An $[n,k,d]_q$ NMDS code can be constructed from an elliptic curve over
$\fq$ having exactly $n$ $\fq$-rational points, for every $k=2,3,\ldots,
n-1$.
%

        \end{theorem}

It should be noted that the proof of Theorem \ref{Th1.1} which appears in
Tsfasman-Vladut book  \cite{Gtsf-vla} depends on
deep algebraic geometry. Here in Section 2 only elementary facts
from algebraic geometry are used to construct certain $[n,k,d]_q$ NMDS codes from
an elliptic curve with $n$ $\fq$-rational points (cf. Theorem \ref{t2.2}).
We will refer to such codes as  \underline{$k$-elliptic codes}.


For every prime power $q$, Theorem \ref{Th1.1} provides  NMDS
codes of length up to $N_q(1)$, where $N_q(1)$ denotes the maximum
number of $\fq$-rational points that an elliptic curve defined
over $\fq$ can have. From work by Waterhouse \cite{DWAT}, we know that  for every $q=p^r$, $p$ prime, $$
N_q(1)= \left \{
\begin{array}{ll}
q+ \lceil 2\sqrt{q} \rceil, & \,\, {\rm for}\,\, p\mid\lceil 2\sqrt{q}\rceil\,\, {\rm and} \,\,
{\rm odd}\,\, r\ge 3,\\
q+\lceil 2\sqrt{q}\rceil +1, & \,\, {\rm otherwise},
\end{array}
\right .
$$
where $\lceil x \rceil$ is the integer part of $x$.

Constructing $[n,k,d]_q$ NMDS codes of length bigger than $N_q(1)$
appears to be hard for $q\ge 17$ and $k\ge 3$ (see \cite{DEBO}).
In Sections 3 and 4 we discuss the related problem whether such
codes can be obtained by extending NMDS $k$-elliptic codes.
In that context the following definition turns out
to be useful.
\begin{definitionn}
{\em An $[n,k,d]_q$ code $\C$ is  \underline{$h$-extendable} if
there exists an $[n+h,k,d+h]_q$ code $\C'$ such that
$\pi_{n,h}(\C')=\C$, where
$\pi_{n,h}:\fq^{n+h}\rightarrow \fq^{n}$,
$\pi_n(a_1,\ldots,a_{n+h})=(a_1,\ldots,a_n)$.
A $1$-extendable code is simply referred to as extendable code.}
\end{definitionn}

With this definition, our main result is stated as follows:

\begin{theorem}\label{main}
Let $q\geq 121$ be an odd prime power. Let $\cE$ be an elliptic
curve defined over $\fq$ whose $j$-invariant $j(\cE)$ is different
from $0$. Then,
\begin{enumerate}
\item for $k=3,6$,
the  $k$-elliptic code associated to $\cE$
is non-extendable;
\item for $k=4$, the  $k$-elliptic code associated to $\cE$
is not $2$-extendable;
\item for $k=5$, the  $k$-elliptic code associated to $\cE$
is not $3$-extendable.
\end{enumerate}
\end{theorem}

\section{Elliptic Codes}
>From now on, $K$ denotes the algebraic closure of the finite field with $q$ elements $\fq$, and
$(X_1, X_2, \ldots, X_k)$ are homogeneous coordinates for $\P^{k-1}(K)$. We also let $X=X_2/X_1$ and
$Y=X_3/X_1$ be the non-homogeneous coordinates for $\P^2(K)$. As usual we identify $(X,Y)\, \in \, K^2$ with the point
$(1,X,Y)\, \in \, \P^2(K)$.

Also, $\cE$ denotes an elliptic plane
curve defined over $\fq$ with affine equation
$$
f(X,Y):=Y^2+a_1XY+a_2Y-X^3-a_3X^2-a_4X-a_5=0\,,
$$
where
$a_i\,\in\,\fq$ for $i=1,\ldots,5$.

Let $n:=\#\cE(\fq)$, the number of $\fq$-rational points of $\cE$.
Then $\cE(\fq)$ consists of $n-1$ affine points, say $P_1, \ldots,
P_{n-1}$, together with its infinite point $P_n=P_{\infty}=(0,0,1)$.

Let $\Sigma=K(x,y)$ be the rational function field of $\cE$, that is the field of
fractions of the domain $K[X,Y]/(f(X,Y))$, where $x=X+(f(X,Y))$ and $y=Y+(f(X,Y))$.
For any point $P\,\in \,\cE$ and for any $\alpha\,\in\,\Sigma$ let $v_P(\alpha)$
denote the order of $\alpha$ in $P$. For $v_P(\alpha)=h>0$, the point $P$ is a zero
of $\alpha$ of multiplicity $h$, and for $v_P(\alpha)=h<0$ the point $P$ is a pole of
$\alpha$ of multiplicity $-h$. By a classical result (see e.g. \cite[Thm.
2.1.50]{Gtsf-vla}), any rational function $\alpha \neq 0$ on an irreducible plane
curve defined over an algebraically closed field has as many zeros as poles, counted
with multiplicity, and $\alpha$ has no zero (and no pole)
if and only if $\alpha$ is constant. As usual, the number of zeros
of $\alpha\in \Sigma$ is indicated by $\ord (\alpha)$. In our case
$\ord (x)=2$, $\ord (y)=3$, $v_{P_\infty}(x)=-2$ and
$v_{P_\infty}(y)=-3$.

For any integer $i>1$, let $$ \psi_i(X,Y):=
\left\{
\begin{array}{ll} Y^s & {\rm
if}\,\, i=3s, \,\,s\ge 1\, ,\\ XY^s & {\rm if}\,\, i=3s+2,
\,\,s\ge 0\, ,\\ X^2Y^s & {\rm if}\,\, i=3s+4,\,\, s\ge 0\,.
\end{array}
\right. $$

Note that $v_{P_\infty}(\psi_i(x,y))=-i$ and hence
ord$(\psi_i(x,y))=i$.

Then, for any $k\,\in\,\{3,4,\ldots, n-1\}$ define the \vspace{.5cm} morphism

$$
\varphi_k:= \left\{
\begin{array}{ccc}
\cE & \rightarrow & \P^{k-1}(K)\\ & & \\
(1,X,Y) & \mapsto & (1,\psi_2(X,Y),\psi_3(X,Y),\ldots, \psi_k(X,Y))
\end{array}
\vspace{.5cm}\right. .$$

Note that $\varphi_k(P_n)=(0,0,\ldots, 0, 1)$.


Let $G_k(\cE)$ be the $(k\times n)$ matrix whose $i^{th}$-column
is the $k$-tuple $\varphi_k(P_i)$ for $i=1,\ldots n$.

\begin{definitionn}\label{neww}
{\em The
subspace of $\fq^{k}$ spanned by the rows of $G_k(\cE)$ is called
the \underline{$k$-elliptic code} associated to $\cE$.}
\end{definitionn}
\noindent \underline{Remark}.$\,\,$
In the notation of \cite{Gtsf-vla}, the $k$-elliptic code associated to $\cE$ is a special Goppa
code, more precisely the code obtained from $(\cE,\cP,D)_L$ by continuation to the point
$P_\infty$ (\cite[p. 271]{Gtsf-vla}), with $\cP=\{P_1,\ldots, P_{n-1}\}$ and  $D=kP_\infty$.

We are in a position to prove the following theorem.

\begin{theorem}\label{t2.2}
For every $k$ with $3\le k \le n-1$, the $k$-elliptic code $\C$
associated to $\cE$ is either an NMDS code or an MDS code of length
$n$ and dimension $k$.
\end{theorem}
\begin{proof} The proof consists of three steps.

\underline{Step 1}. The dimension of $\C$ is equal to $k$ and $d(\C)\ge
n-k$.

For any hyperplane  $\cH$ of $\P^{k-1}(\fq)$, we need to show that
$$ \#(\cH\cap \varphi_k(\cE(\fq))\le k\,. $$

Let $\cH:a_1X_1+a_2X_2+\ldots+a_k X_k=0$. Note that for every $P\, \in\, \cE(\fq)$,
$P\neq P_\infty$, we have that $\varphi_k(P)\, \in \, \cH$ if and only if $P\, \in
\, \cC(\fq)$, where
$\cC$ is the plane curve of equation
$h(X,Y):=a_1+a_2\psi_2(X,Y)+\ldots+a_k\psi_k(X,Y)=0$.

Suppose at first that $a_k\neq 0$, that is $\varphi_k(P_\infty)
\notin \cH$. Then $\#(\cH\cap \varphi_k(\cE(\fq))$ is equal to
the number of affine points in $\cC(\fq)\cap \cE(\fq)$, and hence $\#(\cH\cap
\varphi_k(\cE(\fq))\leq \ord (h(x,y))$. Note that $h\neq 0$,
otherwise $\cE$ would be a component of $\cC$. But this is
impossible, since $h(X,Y)$ has degree in $X$ at most $2$. Then
$v_{P_\infty}(h)\ge -k$, hence $\ord (h)\le k$ and the assertion
follows.

Now, let $a_k=0$. Then we have $\varphi_k(P_\infty)\,\in\, \cH$, whence $\#(\cH\cap
\varphi_k(\cE(\fq))\le 1+\ord (h)$. Again, the assertion follows since
$v_{P_\infty}(h)\ge -(k-1)$ yields $\ord (h)\le k-1$.

\underline{Step 2}. The dimension of $\C^\perp$ is equal to $n-k$ and
$d(\C^\perp)\ge k$.

We need to prove that any $k-1$ points in $\varphi_k(\cE(\fq))$
are linearly independent. Suppose on the contrary that there
exists a set $\cB$ of $k-1$ points in $\varphi_k(\cE(\fq))$
contained in two distinct hyperplanes of $\P^{k-1}(\fq)$, say
$\cH_1:a_1X_1+a_2X_2+\ldots+a_k X_k=0$ and
$\cH_2:b_1X_1+b_2X_2+\ldots+b_k X_k=0$, and consider the rational
functions $h_1:=a_1+a_2\psi_2(x,y)+\ldots+a_k\psi_k(x,y)$ and
$h_2:=b_1+b_2\psi_2(x,y)+\ldots+b_k\psi_k(x,y)$.

If $(0,0,\ldots, 1)\,\notin\, \cB$, then $h_1$ and $h_2$ have at
least $k-1$ common zeros. Moreover, since both $h_1$ and $h_2$ have
order at most $k$, the rational function $h_1/h_2$ has either no or
just one zero. In the former case $h_1/h_2 $ is constant, whence
$\cH_1=\cH_2$, a contradiction. In the latter case, $\ord
(h_1/h_2)=1$, and therefore $\cE$ is isomorphic to $\P^1(K)$, which
is impossible.

Suppose now that $(0,0,\ldots, 1)\,\in\, \cB$. Therefore
$a_k=b_k=0$, hence $\ord (h_1)$ and $\ord (h_2)$ are both less
than or equal to $k-1$, and $h_1$ and $h_2$ have at least $k-2$
zeros in common. This yields $\ord (h_1/h_2)\,\in\,\{0,1\}$ and we
get the same contradiction as above.

\underline{Step 3}. $\C$ is NMDS or MDS.

Step 1 yields that $\C$ is AMDS or MDS. By Step 2 we have $s(\C^\perp)\le
1$, and hence the theorem is proved.
\end{proof}

\noindent \underline{Remark}.$\,\,$
We point out that apart from a few possibilities the
$k$-elliptic code in Theorem \ref{t2.2} is an NMDS code. This is
indeed the case as soon as $\cE$ has $n\geq 5$ $\fq$-rational
points, but a counterexample is known to exist for $n=4$, see
\cite[Thm 3.2.19]{Gtsf-vla}. Here we give an elementary proof under the weaker
hypothesis $n\geq 12$. With same notation as in the proof of
Theorem \ref{t2.2}, we have to prove
$$ \#(\cH\cap
\varphi_k(\cE(\fq))= k\,,
$$ for some hyperplane $\cH$ of $\P^{k-1}(\fq)$. Let $m:=\lceil
\frac{k+1}{3}
\rceil$. We begin by noting that every $h(X,Y)\in
\fq[X,Y]$ of degree $m$ satisfies $$
h(X,Y)-(a_1+a_2\psi_2(X,Y)+\ldots+a_{3m}\psi_{3m}(X,Y))=g(X,Y)f(X,Y) $$ for certain
$a_1,\ldots,a_{3m}\in \fq,\ g\,\in\,K[X,Y].$

Now, take an $\fq$-rational plane curve $\cX$ of order $m$ such
that (i) $\cA:=\cX\cap \cE$ consists of $3m$ $\fq$-rational points
of $\cE$, (ii) $P_\infty\,
\notin\,
\cA$ for $k\equiv 1\pmod 3$ and $P_\infty\, \in\, \cA$ for $k\equiv
-1\pmod 3$. It should be noted that our assumption $n\geq 12$ is
used at this point for the case $m=2$. If $\cX$ has equation
$h(X,Y)=0$ and the coefficients $a_i$ are defined as before, then
the curve of equation
$a_1+a_2\psi_2(X,Y)+\ldots+a_{3m}\psi_k(X,Y)=0$ passes through all
points in $\cA$. Note that the equation
$\cH:a_1X_1+a_2X_2+\ldots+a_{3m}X_{3m}=0 $ defines a hyperplane
$\cH$ for every $k$, since for $k=3m-1$ $P_\infty \, \in \, \cA$
yields $a_{3m}=0$. Then $\cH$ meets $\varphi_k(\cE(\fq))$ in
exactly $k$ points.



\section{Plane elliptic curves and intersections with lines}

The proof of Theorem \ref{main} depends on some results on the
number of $\fq$-rational lines through a given point $P$ which meet
an elliptic cubic curve in exactly three $\fq$-rational points. The
aim of this section is to state and prove such results.

We limit ourselves to the odd order case, that is the underlying
projective plane $\P^2(\fq)$ is assumed to be of odd order $q$.
Then a canonical form for an elliptic cubic curve $\cE$ of
$\P^2(\fq)$ is $Y^2=X^3+aX^2+bX+c$, with $a,b,c\,\in\, \fq$ (see e.g. \cite[p. 46]{Sil}).

We begin with the following lemma.
\begin{lemma}\label{tan}
For every point $P \,\in\,\P^2(\fq)$ not on $\cE$,
\begin{itemize}
\item[\rm (i)] there exist at most $6$ tangents of $\cE$ passing through
$P$;
\item[\rm (ii)] if $P$ is affine, then at least one non-vertical
 line through $P$ is tangent of $\cE$.
\end{itemize}
\end{lemma}
\begin{proof}

The assertion (i) is a classical result in zero characteristic, and
it holds true in positive characteristic $p>3$. So, we may assume
that $p=3$. Now, if the assertion
is false, 
then more
than $6$ tangents to $\cE$ pass through $P$, and hence more than
$6$ points of $\cE$ belong to the polar quadric $\cC$ of $P$ with
respect to $\cE$ (see \cite[Lemma 11.4]{Ch}). Since $\cE$ is
irreducible, B\'ezout Theorem yields that $\cC$ is actually
indeterminate, and hence a line of nuclei of $\cE$ contains $P$
according to \cite[Thm. 11.20(iv)]{Ch}. A straightforward
computation shows that then
$a=b=0$. But this contradicts the non-singularity of $\cE$.

(ii) It is straightforward to check that the intersection between
$\cE$ and the polar quadric of $P=(x_0,y_0)$ with respect to $\cE$
does not entirely consist of points on the line $X=x_0$.
\end{proof}

Let $j(\cE)$ denote the $j$-invariant of the elliptic curve $\cE$.
We start with the case $j(\cE)\neq 0$. The following lemma is an
extension of a result by Hirschfeld and Voloch (\cite[Thm. 5.1]{DHVO}).

\begin{lemma}\label{chiave}
Let $q\geq 121$, and $j(\cE)\neq 0$. Then seven or more lines
through a given $\fq$-rational point $P$ outside $\cE $ intersect
$\cE$ in $3$ distinct $\fq$-rational points.
\end{lemma}
\begin{proof}

Assume at first that $P$ is an affine point, and put $P=(P_x,P_y)$.
Define the rational function $F(X,Y,Z)$ by


$$
-Z^2-Z\big(
a+X-(\frac{Y-P_y}{X-P_x})^2\big)-
\big( X^2+aX+b-2P_y(\frac{Y-P_y}{X-P_x})-\frac{(Y-P_y)^2}{X-P_x}
\big)
$$
Let $Q=(Q_x,Q_y)$ be an $\fq$-rational affine point of $\cE$ such
that $Q_x\neq P_x$. The line through $P$ and $Q$ intersects
$\cE$ in two more (not necessarily distinct) points, say $A$ and
$B$. Then the $X$-coordinates of $A$ and $B$ are roots of the
polynomial $F(Q_x,Q_y,Z)$. In fact, this follows from
$$
F(Q_x,Q_y,Z)=
\frac{1}{Z-Q_x}\big( \big( \frac{Q_y-P_y}{Q_x-P_x}(Z-P_x)+ P_y \big)^2
- Z^3-aZ^2-bZ-c \big)\,.
$$
Next we prove that quadratic polynomial ${\tilde F}(Z)=F(x,y,Z)$ is
irreducible in $\Sigma [Z]$. To do this we may suppose that
$F(x,y,Z)=g(x,y)(Z-h_1(x,y))(Z-h_2(x,y))$, with
$g,h_1,h_1\,\in\,\Sigma$. For $i=1,2$, define the rational \vspace{.5cm} maps
$$
\Phi_i:= \left\{
\begin{array}{ccc}
\cE & \rightarrow & \cE\\ & & \\
(1,X,Y) & \mapsto &
\big(
1,h_i(X,Y),\frac{Y-P_y}{X-P_x}(h_i(X,Y)-P_x) +P_y\big)
\end{array}
\vspace{.5cm}\right. .$$
By
definition of $F$, if $Q=(Q_x,Q_y)\in \cE$ with $Q_x\neq P_x$, then
$\Phi_i(Q)$ belongs to both $\cE$ and the line through $Q$ and $P$.
Moreover, if $\Phi_i$ fixes a point on a non-vertical line
through $P$ then such a line is a tangent of $\cE$. By Lemma \ref{tan}(i),
we have then that $\Phi_i$ has order greater than $4$ or equal to
$3$. Finally, let $l$ be a non-vertical tangent of $\cE$ through
$P$ (such a line exists by Lemma \ref{tan}(ii)). Then, either
$\Phi_1$ or $\Phi_2$ fixes a point in $l\cap \cE$, and therefore
the irreducibility of $F(x,y,Z)$ over $\Sigma(Z)$ follows from
Corollary 4.7 in \cite{Ahar}. Now, we may define the algebraic curve
$\cE'$ as the curve in $\P^3(K)$ whose rational function field is $\Sigma(z)$, $z$ being a root of ${\tilde F}$.
Note that the projection $\pi:\cE'\rightarrow \cE\, , \,\,
\pi(X,Y,Z)=(X,Y)$ is a rational map of degree two.

Suppose that $R=(1,x_1,y_1,z_1)$, $x_1\neq P_x$, is an
$\fq$-rational point of $\cE'$ which is not a ramification point of
$\pi$. Let $\pi^{-1}(\pi(R))=\{R,R'\}$, with $R'=(1,x_1,y_1,z_2)$.
Then $(x_1,y_1)\,\in \, \cE$ and $F(x_1,y_1,z_1)=
F(x_1,y_1,z_2)$ $=$ $ 0$; this means that the line through $P$ and
$(x_1,y_1)$ intersects $\cE$ in three distinct $\fq$-rational
points. Then Lemma \ref{chiave} for an affine point $P$ follows
from the following assertion: The curve $\cE'$ has at least $14$
affine $\fq$-rational non-ramification points $(1,x_1,y_1,z_1)$ such
that $x_1\neq P_x$. To prove it, we note at first that a
ramification point for $\pi$ is a point $(1,x_1,y_1,z_1)$ such that
the line through $P$ and $(x_1,y_1)$ is a tangent to $\cE$. By
Lemma \ref{tan}(i), we may have at most $6$ ramification points.

By Hurwitz Theorem (\cite[Thm. 2.2.36]{Gtsf-vla}) we have that
the genus $g$ of $\cE'$ satisfies $2g-2\le 6$, and hence $g\le 4$.
Let $N$ denote the number of $\fq$-rational points of $\cE'$. By
Hasse-Weil Theorem (\cite[p. 177]{Gtsf-vla}) we have $N\ge
q+1-8\sqrt{q} $, hence $N\ge 34$ from our hypothesis $q\ge 121$.
Then the assertion follows, since $\deg (\cE')=6$ yields that at most $12$ points of
$\cE'$ are in the union of the plane at infinity and the plane of
equation $X=P_x$.

Now assume that $P$ is an infinite point, and put $P=(0,1,m)$. The
proof is similar to the proof given for $P$ affine. Here we define
$$
F_1(x,y,Z):=
\frac{1}{Z-x}\big( \big( m(Z-x)+y \big)^2
- Z^3-aZ^2-bZ-c \big)
$$
instead of $F$. We remark that Lemma \ref{tan}(ii) may not hold for
$P$, since it may happen that the only tangent line through $P$ is
the line at infinity. However, when this occurs, the irreducibility
of ${\tilde F}_1$ still follows from Corollary 4.7 in \cite{Ahar},
since
both $\Phi_1$ and $\Phi_2$  fix the point $(0,0,1)$.
\end{proof}

For $j(\cE)=0$ a result follows from
\cite[Thm 5.2]{FFTA}.

\begin{lemma}\label{mioo}
Let $q=p^r$, $p>3$, $q>9887$. Suppose that $j(\cE)=0$ and that $\cE$ has an even number of
$\fq$-rational points. If $r$ is even or $p\equiv 1 \pmod 3$,
then seven or more lines through a given $\fq$-rational point outside $\cE$ intersect $\cE$ in $3$
distinct $\fq$-rational points.
\end{lemma}

\section{Proof of the  Theorem \ref{main}}
We keep our notation and terminology used in Section 3. Our
approach is based on a strong relationship between $k$-elliptic
codes and certain point-sets in $\P^{k-1}(\fq)$ characterized by
purely combinatorial properties. According to \cite{Chs}, an
\underline{$(n;k,k-2)$-set} in $\P^{k-1}(\fq)$ is defined as a set
consisting of $n$ points no $k+1$ of which lie on the same
hyperplane of $\P^{k-1}(\fq)$. An $(n;k,k-2)$-set in
$\P^{k-1}(\fq)$ is \underline{complete} if it is maximal with respect to
set-theoretical inclusion. From the proof of Theorem \ref{t2.2},
the points of $\varphi_k(\cE(\fq))$ form an $(n;k,k-2)$-set in
$\P^{k-1}(\fq)$.
\begin{lemma}\label{track}
A $k$-elliptic code $\C$ is not-extendable if and only if the
corresponding  $\varphi_k(\cE(\fq))$ is a complete $(n;k,k-2)$-set
in $\P^{k-1}(\fq)$.
\end{lemma}
\begin{proof}
We have to prove that $\C$ is extendable if and only if there exists a point
$P$ in $\P^{k-1}(\fq)\setminus
\varphi_k(\cE(\fq))$ such that no hyperplane through $P$ intersects
$\varphi_k(\cE(\fq))$ in  $k$ points.

Fix a generator matrix for $\C$, say $G_k(\cE)$, and suppose that
no hyperplane through $P\,\in\,\P^{k-1}(\fq)\setminus
\varphi_k(\cE(\fq))$ intersects $\varphi_k(\cE(\fq))$ in $k$
points. Let $G_k(\cE)'$ be the matrix obtained from $G_k(\cE)$ by
adding an extra-column whose entries are the homogeneous
coordinates of $P$. Then the subspace $\C'$ of $\fq^k$ spanned by
the rows of $G_k(\cE)'$ is a $[n+1,k,n-k+1]_q$ code with
$\pi_{n,1}(\C')=\C$.

On the other hand, let $\C'$ be an $[n+1,k,n-k+1]_q$ code with
$\pi_{n,1}(\C')=\C$. Let $R_1=(r_{11}, \ldots, r_{1(n+1)}),\ldots,
R_k=(r_{k1},\ldots ,r_{k(n+1)})$ be an $\fq$-base of $\C'$ such
that $\pi_{n,1}(R_i)$ is the $i$-th row of $G_k(\cE)$. Then no
hyperplane through the point $P=(r_{1(n+1)}, \ldots, r_{k(n+1)})$
intersects $\varphi_k(\cE(\fq))$ in $k$ points.
\end{proof}
Arguing as in Lemma \ref{track}, a more general result can
actually be proved.
\begin{corollary}\label{ctrack}
The $k$-elliptic code $C$ of length $n$ is not $h$-extendable if
the corresponding $(n;k,k-2)$-set $\varphi_k(\cE(\fq))$ is either
complete or can be completed by at most $h-1$ points.
\end{corollary}

We begin the proof of Theorem \ref{main} by noting that the
hypothesis $q\geq 121$ together with the Hasse-Weil theorem ensures
the existence of at least seven $\fq$-rational points on $\cE$. This
shows that $k$-elliptic codes with $k\leq 6$ certainly arise from
$\cE$.

According to Corollary \ref{ctrack}, Theorem \ref{main} will be
proved once we have shown that the $(n;k,k-2)$-set $\varphi_k(\cE(\fq))$ is
either complete or it can be completed by adding at most $h-1$
points where
$$
h:= \left\{
\begin{array}{ll}
1 & {\rm for }\,\, k=3,6\, ;\\ 2 & {\rm for }\,\, k=4\, ;\\ 3 &
{\rm for}\,\, k=5\,.
\end{array}
\right.  $$
Lemma \ref{chiave} allows us to choose a frame in
$\P^2(\fq)$ satisfying the following conditions:
\begin{itemize}
\item the line of equation $X=0$ meets $\cE$ in two affine $\fq$-rational points,
both distinct from $(0,0)$;
\item both lines $Y=0$ and $X=Y$ meet $\cE$ in three affine $\fq$-rational points.
\end{itemize}

We distinguish several cases according to the value of $k$.

\underline{Case} $k=3$. \\ By Lemma \ref{chiave}, $\varphi_3(\cE(\fq))$ is
complete.

\underline{Case} $k=4$. \\ Let $\varphi_4(\cE(\fq))$ be incomplete,
and choose a point $Q=(Q_1,Q_2,Q_3,Q_4)$ in $\P^{3}(\fq)$ that
can be added to $\varphi_4(\cE(\fq))$.
We show that such a point $Q$ lies on the line through the
fundamental points $(0,0,1,0)$ and $(0,0,0,1)$. In fact,
for $(Q_1,Q_2,Q_3)\neq(0,0,1)$, Lemma \ref{chiave} implies the
existence of a line $l:a+bX+cY=0$ through $P=(Q_1,Q_2,Q_3)$ that
meets $\cE$ in three distinct $\fq$-rational affine points. Then
the plane of equation $aX_1+bX_2+cX_3+0X_4$ passes through $Q$ and
meets $\varphi_4(\cE)$ in $4$ distinct $\fq$-rational points,
more precisely the points in $\{\varphi_4(l\cap
\cE(\fq)),(0,0,0,1)\}$. But this is impossible since $Q$ is
assumed to be a point that can be added to $\varphi_4(\cE(\fq))$.
This contradiction proves the assertion. Now, to prove
Theorem \ref{main} for $k=4$, we have to check that
$\varphi_4(\cE(\fq))\cup
\{Q\}$ is complete, that is no further point $Q'=(0,0,1,\beta)$,
$\beta \,\in\, \fq$, can be added to $\varphi_4(\cE(\fq))\cup
\{Q\}$. But this follows immediately from the fact that the plane $X_2=0$ passes
through $Q',\,Q$ and three distinct points in $\varphi_4(\cE(\fq))$, which are those in
$\{\varphi_4(\{X=0\}\cap \cE(\fq)),(0,0,0,1)\}$.

\underline{Case} $k=5$. Let
$Q=(Q_1,Q_2,Q_3,Q_4,Q_5)\,\in\,\P^4(\fq)\setminus \varphi_5(\cE(\fq))$.
We need the following technical lemma.

\begin{lemma}
\label{k=5}
If $Q$ can be added to
$\varphi_5(\cE(\fq))$, then $Q_5Q_2\neq 0$, $Q_4=0$ and
$(1,0,Q_5/Q_2)\,\in\, \cE$.
\end{lemma}
\begin{proof}
If $Q_5=0$, then the hyperplane $X_5=0$ meets $\varphi_5(\cE)$ in
$5$ distinct $\fq$-rational points, which are those in
$\{\varphi_5(\{XY=0\}\cap \cE(\fq))\}$.

For $Q_5\neq 0$, $Q_2=0$, $Q_4=0$, Lemma \ref{chiave} ensures the
existence a line $l$ through $P=(0,0,1)$ which is different from
$X=0$ and meets $\cE$ in two more distinct $\fq$-rational affine
points. If $l$ has equation $X+\alpha=0$, then the hyperplane in
$\P^4(\fq)$ of equation $\alpha X_2+X_4=0$ passes through $Q$ and
meets $\varphi_5(\cE)$ in $5$ distinct $\fq$-rational points, which are those in
$\{\varphi_5(\{X(X+\alpha)=0\}\cap \cE(\fq)),
(0,0,0,0,1)\}$.

Similarly, for $Q_5\neq 0$, $Q_2=0$, $Q_4\neq 0$: A line $l$ through
$P=(0,Q_4/Q_5,1)$ meets $\cE$ in three distinct $\fq$-rational affine
points not lying on $X=0$. If $l:\alpha(X-Q_4/Q_5 Y)+\beta =0$,
then the hyperplane of equation $\beta X_2+ \alpha X_4-\alpha
Q_4/Q_5 X_5=0$ passes through $Q$ and meets $\varphi_5(\cE)$ in
$5$ distinct $\fq$-rational points. Also, for $Q_5\neq 0$, $Q_2\neq 0$,
$Q_4\neq 0$: A line of equation $\alpha(X-Q_4/Q_2)+\beta
(Y-Q_5/Q_2)=0$ meets $\cE(\fq)$ in three distinct $\fq$-rational affine
points not lying on $X=0$, and
%
%
the hyperplane $\alpha (X_4-Q_4/Q_2 X_2)+ \beta (X_5- Q_5/Q_2
X_2)=0$ passes through $Q$ and meets $\varphi_5(\cE)$ in $5$
$\fq$-rational points. Finally for $Q_5\neq 0$, $Q_2\neq 0$, $Q_4=
0$, $(1,0,Q_5/Q_2)\,\notin\,\cE$: A line of equation $\alpha
X+\beta (Y-Q_5/Q_2)=0$ meets $\cE(\fq)$ in three $\fq$-rational
affine points not lying on $X=0$, and
%
the hyperplane $\alpha X_4+ \beta (X_5- Q_5/Q_2 X_2)=0$ passes
through $Q$ and meets $\varphi_5(\cE)$ in $5$ distinct
$\fq$-rational points. This completes the proof of Lemma \ref{k=5}.
\end{proof}

To settle the case $k=5$ suppose that $Q$ can be added to
$\varphi_5(\cE(\fq))$. Let $\{X=0\}\cap
\cE=\{(0,0,1), (1,0,
\lambda), (1,0, \mu)\}$, and assume $\lambda = Q_5/Q_2$.

Note that no point $Q'=(Q_1',Q_2',Q_3',0,Q_5')$, with $Q_2'Q_5'\neq
0$ and such that  $Q_5'/Q_2'=\lambda$ can be added to $\varphi_5(\cE(\fq))
\cup \{Q\}$. Lemma \ref{chiave} ensures the existence of a line $l$ through
$P=(1,0,\lambda)$ that meets $\cE$ in three distinct $\fq$-rational affine
points, two of which not lying on $X=0$. If $l:\alpha X+\beta
(Y-\lambda)=0$, then the hyperplane of equation $\alpha X_4+ \beta
(X_5- \lambda X_2)=0$ passes through $Q'$ and meets
$\varphi_5(\cE(\fq)) \cup \{Q\}$ in $5$ distinct points.

This shows that if a point $Q'$ can be added to
$\varphi_5(\cE(\fq))\cup
\{Q\}$ then $Q'=(Q_1',1,Q_3', 0, \mu)$. Finally, a straightforward argument
shows that $\varphi_5(\cE(\fq))\cup
\{Q,Q'\}$ is complete.

\underline{Case} $k=6$.

Given any point
$Q=(Q_1,Q_2,Q_3,Q_4,Q_5,Q_6)\,\in\,\P^5(\fq)\setminus
\varphi_6(\cE)$, we have to find a hyperplane $\cH$ of
$\P^5(\fq)$ through $Q$ that meets $\varphi_6(\cE)$ in 6 distinct
$\fq$-rational points. To do this, we distinguish a number of
cases, even if we use the same kind of argument depending on Lemma
\ref{chiave}.

1) $Q_5=0$.  The hyperplane $X_5=0$ passes through $Q$ and meets $\varphi_6(\cE)$ in $6$
distinct $\fq$-rational points, which are those in
$\{\varphi_6(\{XY=0\}\cap \cE(\fq)), (0,0,0,0,0,1)\}$.

2) $Q_5=1,\,Q_4=Q_6=0,\,Q_2\neq Q_3$. Let $l$ be a line
through $P=(1, \frac{1}{Q_3-Q_2},-\frac{1}{Q_3-Q_2})$
meeting $\cE$ in three distinct $\fq$-rational points outside the line
$X=Y$. If $l$ has equation $\alpha (1+(Q_3-Q_2)Y)+\beta (X+Y)=0$,
then the hyperplane $\alpha (X_2-X_3)+\beta X_4+\alpha
(Q_3-Q_2)X_5+(-\beta
-\alpha (Q_3-Q_2))X_6=0$ passes through $Q$ and meets
$\varphi_6(\cE)$ in $6$ distinct $\fq$-rational points, more precisely the points in
$\{\varphi_6((\{X-Y=0\}\cup l) \cap \cE(\fq)) \}$.

3) $Q_5=1,\,Q_4=Q_6=0,\,Q_2=Q_3$.
A line of equation $\alpha+\beta (X+Y)=0$
meets $\cE$ in three distinct $\fq$-rational points outside the line $X=Y$.
Then the hyperplane of equation $\alpha (X_2-X_3)+\beta
X_4-\beta X_6=0$ passes through $Q$ and meets $\varphi_6(\cE)$ in
$6$ distinct $\fq$-rational points.

4) $Q_5=1,\,Q_6\neq 0, Q_3=0$.
A line of equation $\alpha +\beta (X-Y/Q_6)=0$
meets $\cE$ in three distinct $\fq$-rational points outside the line
$Y=0$. Then the hyperplane of equation $\alpha X_3+\beta X_5-\beta/Q_6 X_6=0$
passes through $Q$ and meets $\varphi_6(\cE)$ in $6$
distinct $\fq$-rational points.

5) $Q_5=1,\,Q_6\neq 0, Q_3\neq 0$.
A line of equation $\alpha (X- 1/Q_3) +\beta
(Y-Q_6/Q_3)=0$ meets $\cE$ in three distinct $\fq$-rational points outside
the line $Y=0$, and  the hyperplane $\alpha (X_5- X_3/Q_3)
+\beta (X_6-Q_6/Q_3X_3)=0$ passes through $Q$ and meets
$\varphi_6(\cE)$ in $6$ distinct $\fq$-rational points.

6) $Q_5=1,\,Q_4\neq 0, Q_2=0$.
A line of equation $\alpha (X-Q_4Y)+\beta=0$ meets $\cE$
in three distinct $\fq$-rational points not lying on the line $X=0$. Then the
hyperplane $\alpha (X_4-Q_4X_5)+\beta X_2=0$ passes through $Q$
and meets $\varphi_6(\cE)$ in $6$ distinct $\fq$-rational points, which are those
in $\{\varphi_6((\{X=0\}\cup l) \cap \cE(\fq)),
(0,0,0,0,0,1)\}$.

7) $Q_5=1,\,Q_4\neq 0, Q_2\neq 0$.
A line of equation $\alpha
(X-Q_4/Q_2)+\beta (Y-1/Q_2)=0$ meets $\cE$ in three distinct
$\fq$-rational points outside the line $X=0$, and the hyperplane
$\alpha (X_4-Q_4/Q_2X_2)+\beta (X_5-X_2/Q_2)=0$ passes through
$Q$ and meets $\varphi_6(\cE)$ in $6$ distinct $\fq$-rational points.

As a consequence of Lemma \ref{mioo}, an analogous to Theorem
\ref{main} can be proved for some cubics $\cE$ with $j(\cE)=0$.

\begin{theorem}
Let $q=p^r$, $p>3$, $q>9887$. Let $\cE$ be an elliptic curve defined over $\fq$, with $j(\cE)=0$ and
having an even number of $\fq$-rational points.
If $r$ is even or $p\equiv 1 \pmod 3$, then
\begin{enumerate}
\item for $k=3,6$,
the $k$-elliptic code associated to $\cE$ is non-extendable;
\item for $k=4$, the $k$-elliptic code associated to $\cE$
is not $2$-extendable;
\item for $k=5$, the $k$-elliptic code associated to $\cE$
is not $3$-extendable.
\end{enumerate}
\end{theorem}
\noindent \underline{Remark}.$\,\,$
Our method still works for $k>6$ even if some modification is needed. However, the result is not so
sharp as for $k\le 6$ since it only ensures non-$h$-extendability for $h$ sufficiently bigger than
$k$.

\section*{Acknowledgments}
The author would like to thank Prof. G. Korchm\'aros and Prof. F. Torres for their useful comments.


\begin{thebibliography}{1}
















\bibitem{CHEN} Chen, H., Yau, S.S.-T.: Contribution to Munuera's Problem on the Main Conjecture of Geometric
Hyperelliptic MDS Codes. IEEE Trans. Inform. Theory {\bf 43} 1349--1354, (1997).







\bibitem{DEBO} De Boer, M.A.:
Almost MDS codes. Des. Codes Cryptogr. {\bf 9}, 143--155 (1996).

\bibitem{DDEC} Di Comite, C.:
Su $k$-archi deducibili da cubiche piane.
Atti Accad. Naz. Lincei Rend. {\bf 33}, 429--435 (1962).

\bibitem{DDEC2} Di Comite, C.:
Intorno a certi $(q+9)/2$-archi completi di $S_{2,q}$.
Atti Accad. Naz. Lincei Rend. {\bf 36}, 819--824 (1964).

\bibitem{DODU} S.M. Dodunekov, S.M., Landjev, I.N.:
On near-MDS codes.  J. Geom. {\bf 54}, 30--43 (1995).

\bibitem{DRIE} Driencourt, Y., Michon, J.F.: Remarques sur les codes geometriques.
C.R. Acad. Sci., Paris, Ser. I {\bf 301}, 15--17 (1986).














\bibitem{FUL} Fulton, W.: Algebraic Curves. An introduction to Algebraic
Geometry. New York: W.A. Benjamin 1969.












\bibitem{FFTA} Giulietti, M.: On plane arcs contained in cubic curves. Finite Fields
 Appl. {\bf 8}, 69--90  (2002).











\bibitem{Ahar} Hartshorne, R.: Algebraic Geometry. Grad. Texts in Math. Vol. 52. Berlin
New York: Springer Verlag 1977.





\bibitem{ARS} Hirschfeld, J.W.P.:  Finite Projective Spaces of Three
dimension. Oxford: Oxford University Press 1985.

\bibitem{Ch} Hirschfeld, J.W.P.:  Projective Geometries Over Finite
Fields, 2nd edition. Oxford: Oxford University Press 1998.








\bibitem{Chs} Hirschfeld, J.W.P., Storme, L.: The packing problem in
statistics, coding theory and finite projective planes. J. Statist.
Plann. Inference {\bf 72}, 355--380 (1998).


\bibitem{CHT} Hirschfeld J.W.P., Thas, J.A.: General Galois
Geometries. Oxford: Clarendon Press 1991.

\bibitem{DHVO} Hirschfeld J.W.P., Voloch, J.F.:
The characterization of elliptic curves over finite fields.
J. Austral. Math. Soc. Ser. A {\bf 45}, 275--286 (1988).














\bibitem{Clan} Landjev, I.N.: Linear codes over finite fields and
finite projective geometries. Discrete Math. {\bf 213}, 211--244
(2000).






\bibitem{MAC} MacWilliams, F.J., Sloane, N.J.: The theory of error-correcting
codes. Amsterdam: North-Holland 1977.

\bibitem{MMP} Marcugini, S., Milani, A., Pambianco, F.: Existence and classification of NMDS
codes over $GF(5)$ and $GF(7)$. Proceedings of ACCT2000, the Seventh International Workshop on
Algebraic and Combinatorial Coding Theory, Bulgaria, June 2000, 232--239.


\bibitem{MMP3}  Marcugini, S., Milani, A., Pambianco, F.:
MNDS codes of maximal length over $GF(q)$, $8\le q \le 11$. IEEE Trans. Inform. Theory, submitted.


\bibitem{MUNU1} Munuera, C.: On the main conjecture on geometric MDS
codes. IEEE Trans. Inform. Theory {\bf 38}, 1573--1577 (1992).


\bibitem{MUNU2} Munuera, C.: On MDS elliptic codes. Discrete Math.
{\bf 117} 279--286 (1993).





%




















\bibitem{Schrok} Shokrollahi, M.A.: Minimum Distance of Elliptic Codes
Adv. in Math. {\bf 93} 251--281 (1992).


\bibitem{Sil} Silverman, J.H.:
The Arithmetic of Elliptic Curves. Berlin Heidelbeg New York
Tokio: Springer Verlag 1986.


\bibitem{CSTI} Stichtenoth, H.:
Algebraic Function Fields and Codes.
Berlin Heidelbeg New York
Tokio: Springer Verlag 1993.














\bibitem{DSZ3} Sz\H onyi, T.:
Arcs in cubic curves and 3-independent subsets of abelian groups.
Combinatorics, Eger, Colloq. Math. Soc. J\'anos Bolyai {\bf
52} , 499--508 (1987).








\bibitem{Gtsf-vla}  Tsfasman, M.A., Vladut, S.G.:  Algebraic-Geometric
Codes. Amsterdam: Kluwer 1991.




\bibitem{DVOL} Voloch, J.F.:
On the completeness of certain plane arcs.
European J. Combin. {\bf 8}, 453--456 (1987).

\bibitem{DVO2} Voloch, J.F.:
On the completeness of certain plane arcs II.
European J.  Combin. {\bf 11}, 491--496 (1990).




\bibitem{DWAT} Waterhouse, W.G.:
Abelian varieties over finite fields.
Ann. Sci. \'Ecole Norm. Sup. {\bf 2}, 521--560 (1969).


\bibitem{DZIR} Zirilli, F.:
Su una classe di k-archi di un piano di Galois. Atti Accad.
Naz. Lincei Rend. {\bf 54} , 393--397 (1973).






\end{thebibliography}

%

\end{document}